\documentclass{siamltex}
\usepackage{a4}
\usepackage[dvips]{graphicx}
\usepackage{amsmath}

\renewcommand{\theequation}{\thesection.\arabic{equation}}
\newtheorem{rmrk}{Remark}[section]

\title{Symmetry reductions of a nonlinear option pricing model}

\author{L. A. Bordag\thanks{Halmstad University, Box 823, 301 18 Halmstad, Sweden({\tt Ljudmila.Bordag@ide.hh.se}) }}
\begin{document}
\maketitle
\begin{abstract}
 The studied model was suggested to design a perfect hedging strategy for a
 large trader. In this case the
 implementation of a hedging strategy affects the price of the underlying
 security. The feedback-effect leads to a nonlinear version of the
 Black-Scholes partial differential equation. 
Using the Lie group theory we reduce the
partial differential equation in special cases to ordinary
differential equations. The found Lie group of the model equation gives
 rise to invariant solutions.
Families of exact invariant solutions for special values of parameters  are described.
\end{abstract}



\setlength{\parindent}{0cm}

\pagestyle{myheadings}
\thispagestyle{plain}
\markboth{L. A. BORDAG}{Symmetry reductions of a nonlinear model}

\section{Introduction}

 In a series of works \cite{Frey:perfect}, \cite{Frey:market},
 \cite{Frey:illicuidity} and \cite{FreyPatie} a model for a hedging
 strategy in an illiquid market was suggested.  In the model the
 implementation of a hedging strategy affects the price of the
 underlying security. For a large trader a hedge-cost of the claim
 differs from the price of the option. The feedback-effect leads to a
 nonlinear version of the Black-Scholes partial differential equation,
\begin{eqnarray} \label{intur}
u_t+\frac{\sigma^2 S^2}2\frac{u_{SS}}{(1-\rho S \lambda (S) u_{SS})^2}=0,
\end{eqnarray}
with $S \in [0,\infty), ~~ t \in [ 0,T ]. $ As usual, $S$ denotes here
the price of the underlying asset and $u(S,t)$ denotes the hedge-cost
of the claim with later defined payoff, which is different from the
price of the derivatives product in illiquid markets, $t$ is the time
variable, $\sigma$ defines the volatility of the underlying asset,
$\rho$ is a measure for the feedback-effect of a large trader,
$\lambda (S)$ is chosen in a way to obtain the desired payoff.  The
values of $\rho$ and $\lambda(S)$ might be estimated from the observed
option prices. In dependence on the propositions on the market
different variations of the Black-Scholes formula can accrue like in a
well known model \cite{SchonbucherWilmotta}. Usually the volatility
term in the Black-Scholes formula will be replaced to fit the behavior of the 
price on the market. 
The modeling
process is not finished now and new models can appear. An analytical study 
of these equations may be useful for an easier classification of models created.

 Frey and co-authors  studied  equation
(\ref{intur}) under constrictions and did some numerical simulations.
Our goal is  to investigate this equation using analytical methods.\\
We study the model equation (\ref{intur}) using methods of the Lie
group theory. This method has a long tradition beginning with the work
of S.~Lie \cite{Lie}.  The applications of this method are connected
with an obvious limitation of group-theoretical methods based on
local symmetries because many nonlinear partial differential equations
do not have local symmetries. The modern description of the method and
a large number of applications can be found in \cite{Ovsiannikov},
\cite{Olver}, \cite{Stephani}, \cite{Gaeta}, \cite{Ibragimov}.

In Section \ref{liesec2} we find the Lie algebra and finite equations
for the symmetry group of equation (\ref{intur}). For a special form of the
function $\lambda(S)$ it is  possible to find two functionally independent
invariants of the symmetry group.\\
Using the symmetry group and its invariants we reduce the
partial differential equation (\ref{intur}) in special cases to ordinary
differential equations in Section \ref{redsec}. 
We study singular points of the
reduced equations in Section \ref{k1expl}
and describe the behavior of invariant solutions. 
For a fixed set of parameters the complete set of exact invariant solutions
is given.

\section{Lie group symmetries}\label{liesec2}

Let us introduce a two-dimensional space $X$ of
independent variables $(S,t) \in X$ and a one-dimensional space of 
dependent variables $u \in U.$ We consider the space $U_{(1)}$ of the
first derivatives of the variable $u$ on $S$ and $t$, i.e., $(u_S,u_t)
\in U_{(1)}$ and analogously we introduce the space $U_{(2)}$ of the
second order derivatives $(u_{SS}, u_{St}, u_{tt}) \in U_{(2)}.$ We
denote by $M=X \times U$ a base space which is a Cartesian product of
pairs $(x,u)$ with $x=(S,t) \in X,~~u\in U$.  The studied differential
equation (\ref{intur}) is of the second order and to represent this
equation as an algebraic equation we introduce a
second order jet bundle $ M^{(2)}$ of the base space $M$. 
This space has the form
\begin{equation}
M^{(2)} = X \times U \times U_{(1)} \times U_{(2)} \label{jet}
\end{equation}
and possesses a natural contact structure.
We label the coordinates in the space $M^{(2)}$ by $w=(S,t,u,u_S,u_t,u_{SS}, u_{St},
u_{tt}) \in M^{(2)}.$

In the space $M^{(2)}$ equation (\ref{intur}) is equivalent to the relation
\begin{equation}\label{algeq}
\Delta(w)=0 , ~~ w \in M^{(2)},
\end{equation}
where we denote by $\Delta$
  the following function
\begin{equation}
\Delta(S,t,u,u_S,u_t,u_{SS}, u_{St}, u_{tt})=u_t+\frac{\sigma^2
  S^2}2\frac{u_{SS}}{(1-\rho S\lambda(S) u_{SS})^2}. \label{delta}
\end{equation}
We identify the algebraic equation (\ref{algeq}) with its solution
manifold $L_{\Delta}$ defined by
\begin{equation}
L_{\Delta}=\{ w \in M^{(2)} | \Delta(w)=0 \} \subset M^{(2)}.
\end{equation}

Let us consider an action of a Lie-point group on our differential 
equation
and its solutions. 
We define a symmetry group $G_\Delta$ of equation (\ref{algeq}) by
\begin{equation}
G_\Delta=\{ g \in \rm{Diff}( M^{(2)})|~~ g: ~~L_{\Delta} \to L_{\Delta}\},
\end{equation}
consequently we are interested in a subgroup of $\rm{Diff}( M^{(2)})$ which is
compatible with the structure of $L_{\Delta}$.

As usual we first find the corresponding symmetry Lie algebra ${\mathcal
Diff}_{\Delta} ( M^{(2)}) \subset {\mathcal
Diff} ( M^{(2)})$ and then use the main Lie theorem to obtain
$G_\Delta$ and its invariants.

We denote 
 an element of a Lie-point vector field on $M$
by
\begin{equation}\label{vekv}
V= \xi (S,t,u) \frac{\partial}{\partial S} + \tau (S,t,u)
\frac{\partial}{\partial t} + \phi(S,t,u) \frac{\partial}{\partial u},
\end{equation}
where $\xi (S,t,u)$,$\tau (S,t,u)$ and $\phi(S,t,u)$ are smooth
functions of their arguments, $V \in {\mathcal Diff} (M) $. 

If the infinitesimal generators of  $g \in G_{\Delta} $
exist then they have the structure of the type (\ref{vekv}) and form an algebra
${\mathcal Diff}_{\Delta} ( M).$

 A Lie group of transformations acting on the base space $M$ induce 
as well the transformations on $M^{(2)}.$ 

The
corresponding algebra ${\mathcal Diff}_{\Delta} ( M^{(2)})$ will be
composed of vectors
\begin{eqnarray}
pr^{(2)} V &=& \xi (S,t,u) \frac{\partial}{\partial S} + \tau (S,t,u)
\frac{\partial}{\partial t} + \phi(S,t,u) \frac{\partial}{\partial u} \nonumber \\
&+&\phi^S(S,t,u) \frac{\partial}{\partial u_S}+\phi^t(S,t,u) \frac{\partial}{\partial u_t}  \label{prol2} \\
&+&\phi^{SS}(S,t,u) \frac{\partial}{\partial u_{SS}}+\phi^{St}(S,t,u) \frac{\partial}{\partial u_{St}}
+\phi^{tt}(S,t,u) \frac{\partial}{\partial u_{tt}} ,\nonumber
\end{eqnarray}
which are the second prolongation of vectors $V$.  Here the smooth
functions $\phi^S(S,t,u)$, $\phi^t(S,t,u)$, $\phi^{SS}(S,t,u)$,
$\phi^{St}(S,t,u)$ and $\phi^{tt}(S,t,u)$ are uniquely defined by the
functions $ \xi (S,t,u), \tau (S,t,u)$ and $\phi(S,t,u)$ using the
prolongation procedure (see \cite{Ovsiannikov}, \cite{Olver},
\cite{Stephani}, \cite{Gaeta}, \cite{Ibragimov}).

\begin{theorem}
The differential equation (\ref{intur}) with an arbitrary function
  $\lambda(S)$ possesses a trivial three dimensional Lie algebra $ Diff_{\Delta} (M)$ spanned by
  generators $$ V_1 =  \frac{\partial}{\partial t}, ~~
V_2 =  S\frac{\partial}{\partial u},~~
V_3 = \frac{\partial}{\partial u}.$$ 
Only for the special form of the function 
$\lambda(S) \equiv \omega S^k,$ where $ \omega, k\in { R}$ equation (\ref{intur})
  admits a nontrivial four dimensional Lie algebra spanned by generators
$$ V_1 =  \frac{\partial}{\partial t}, ~~
V_2 =  S\frac{\partial}{\partial u},~~
V_3 = \frac{\partial}{\partial u}, 
V_4 = S \frac{\partial}{\partial S}+ (1-k) u \frac{\partial}{\partial u}.$$ 
\end{theorem}

\begin{proof}
The symmetry algebra ${\mathcal Diff}_{\Delta} (
M^{(2)})$ of the second order differential equation (\ref{algeq}) can be
found as a solution of the defining equations
\begin{equation}
pr^{(2)} V(\Delta)=0 ~(mod(\Delta =0)), \label{algebradef}
\end{equation}
i.e., the equation (\ref{algebradef}) should be satisfied on the solution manifold $L_{\Delta}$.

For our calculations we will use the exact form of the coefficients $\phi^t(S,t,u)$ and $\phi^{SS}(S,t,u)$ only.
The coefficient $\phi^t(S,t,u)$ can be defined by the formula
\begin{equation}
\phi^t(S,t,u)= \phi_t+ u_t \phi_u -u_S \xi_t - u_S u_t \xi_u - u_t \tau_t - (u_t)^2 \tau_u, \label{fit}
\end{equation}
and the coefficient $\phi^{SS}(S,t,u)$ by the expression
\begin{eqnarray}
\phi^{SS}(S,t,u)&=&\phi_{SS} + 2 u_S \phi_{S u} + u_{SS} \phi_u \label{fiss} \\
&+&(u_S)^2 \phi_{uu} - 2
 u_{SS} \xi_S - u_S \xi_{SS}- 2 (u_S)^2 \xi_{Su}\nonumber \\
&-&3 u_S u_{SS} \xi_u -(u_S)^3 \xi_{uu} -2 u_{St} \tau_S -u_t \tau_{SS}\nonumber \\
&-&2 u_S u_t \tau_{Su}-(u_t u_{SS}+2 u_S u_{St}) \tau_{u} - (u_S)^2 u_t \tau_{uu} \nonumber ,
\end{eqnarray}
where the subscripts by $\xi, \tau , \phi$ denote corresponding
partial derivatives.  

The first equations of the set (\ref{algebradef}) imply that if $V\in {\mathcal Diff}_{\Delta} ( M)$ then 
\begin{equation}
\xi (S,t,u)= a_1 S,~~\tau (S,t,u)=a_2,~~\phi(S,t,u)=
a_3 S + a_4 + a_5 u, 
\end{equation}
where $a_1,a_2,a_3,a_4,a_5$ are arbitrary constants and $\xi$,$\tau$, $\phi$ are coefficients in the expression (\ref{vekv}).

The remaining equation has a form
\begin{equation}\label{poslmog}
a_1 S \lambda_S(S) -(a_1-a_5)\lambda(S)=0.
\end{equation}
Because this equation should be satisfied for all $S$ identically we obtain for an arbitrary function $\lambda(S)$
\begin{equation}\label{allfal}
a_1=a_5=0,~~\to \xi (S,t,u)=0,~~\tau (S,t,u)=a_2,~~\phi(S,t,u)=a_3 S + a_4.
\end{equation}

Finally, ${\mathcal Diff}_{\Delta}(M)$ admits the following generators
\begin{eqnarray}
V_1 =  \frac{\partial}{\partial t}, ~~
V_2 &=&  S\frac{\partial}{\partial u},~~
V_3 = \frac{\partial}{\partial u}, \label{al3}
\end{eqnarray}
with commutator relations 
\begin{eqnarray}
[V_1,V_2]=[V_1,V_3]=[V_2,V_3]=0.
\end{eqnarray}

If the function $\lambda(S)$ has a special form
\begin{equation}\label{lam}
\lambda(S) \equiv \omega S^k,~~~ \omega, k\in { R}
\end{equation}
then the equation (\ref{poslmog}) on the coefficients of (\ref{vekv}) is less restrictive and we
obtain
\begin{equation}\label{specfal}
 \xi (S,t,u)=a_1 S,~~\tau (S,t,u)=a_2,~~\phi(S,t,u)=(1-k)a_1 u+ a_3 S + a_4.
\end{equation}
Now the symmetry algebra ${\mathcal Diff}_{\Delta} ( M)$
admits four generators
\begin{eqnarray}
V_1 =  \frac{\partial}{\partial t}, ~~
V_2 &=&  S\frac{\partial}{\partial u},~~
V_3 = \frac{\partial}{\partial u}, \label{al4}
V_4 = S \frac{\partial}{\partial S}+ (1-k) u \frac{\partial}{\partial u},
\end{eqnarray}
with commutator relations 
\begin{eqnarray}
~[V_1,V_2]=[V_1,V_3]=[V_1,V_4]=[V_2,V_3]=0,\nonumber\\
~[V_2,V_4]=-k V_2,~~~~[V_3,V_4]=(1-k) V_3. \label{alge}
\end{eqnarray}
\qquad\end{proof} 

{\rmrk In the general case the algebra (\ref{{al4}})
possesses a two dimensional Abelian sub-algebra. For the cases $k=0,1$
the Abelian sub-algebra is three dimensional (\ref{alge}) and we see
later that the corresponding equations (\ref{urav1}) became
autonomous.}
\vspace{3pt}

The symmetry algebra ${\mathcal Diff}_{\Delta} ( M)$ defines by the
main theorem of S. Lie \cite{Lie} the corresponding symmetry
group $G_\Delta$ of the equation (\ref{algeq}).  
To find the closed form of
transformations for the solutions of equation (\ref{intur})
corresponding to this symmetry group we just integrate the system of
ordinary differential equations
\begin{eqnarray}\label{symsys1}
\frac{d{\tilde S}}{d\epsilon}=\xi ({\tilde S},{\tilde t},u), \\
~\frac{d{\tilde t}}{d\epsilon}=\tau ({\tilde S},{\tilde t},{\tilde u}),\\
~\frac{d{\tilde u}}{d\epsilon}=\phi ({\tilde S},{\tilde t},{\tilde u}),
\end{eqnarray}
with initial conditions
\begin{equation}\label{gransys1}
{\tilde S}|_{\epsilon=0}=S,~{\tilde t}|_{\epsilon=0}=t,~{\tilde u}|_{\epsilon=0}=u.
\end{equation}
Here the variables ${\tilde S},{\tilde t}$ and ${\tilde u}$
denote  values $S,t,u$ after a symmetry transformation. The parameter
$\epsilon$ describes a motion along an orbit of the group.

\begin{theorem} 
\label{symteor}
The action of the symmetry group $G_{\Delta}$ of (\ref{intur}) with an arbitrary function
  $\lambda(S)$ is given by (\ref{str})--(\ref{utro}). If the function $\lambda(S)$
has the special form (\ref{lam}) then the symmetry group $G_{\Delta}$ is
represented by (\ref{strsp})-(\ref{utrosp}).
\end{theorem}

\begin{proof}
The solutions of the system of ordinary differential equations
(\ref{symsys1}) with functions $\xi$,$\tau$, $\phi$ defined by (\ref{allfal})
and initial conditions (\ref{gransys1}) have the
form
\begin{eqnarray}
{\tilde S}&=&S,\label{str}\\
{\tilde t}&=&t+ a_2 \epsilon,\label{ttr}\\
{\tilde u}&=&u +a_3 S \epsilon +a_4 \epsilon,~~ \epsilon \in (- \infty,\infty). \label{utro}
\end{eqnarray}
The equations (\ref{str})--(\ref{utro}) are the finite
representation of the symmetry group $G_\Delta$ which corresponds to the symmetry algebra 
 defined by (\ref{al3}) in case of an arbitrary function  $\lambda(S)$.\\
If the function $\lambda(S)$ has a special form given by (\ref{lam})
we obtain a reacher symmetry group. The solution of the system of equations 
(\ref{symsys1}) with the functions 
$\xi$,$\tau$, $\phi$ defined by (\ref{specfal}) and initial conditions
(\ref{gransys1}) have the
form
\begin{eqnarray}
{\tilde S}&=&S e^{a_1 \epsilon},~~ \epsilon \in (- \infty,\infty ),\label{strsp}\\
{\tilde t}&=&t+ a_2 \epsilon,\label{ttrsp}\nonumber\\
{\tilde u}&=&u e^{a_1 (1-k) \epsilon} + \frac{a_3}{a_1 k} S \epsilon e^{a_1 \epsilon} (1-e^{- a_1 k\epsilon} )   \nonumber\\
&+&\frac{a_4}{a_1(1-k)}(e^{a_1 (1-k) \epsilon} -1),~~ k \ne 0,~k\ne 1 
\label{utrnesp}\\
{\tilde u}&=&u e^{a_1\epsilon} + a_3 S \epsilon e^{a_1 \epsilon} +
\frac{a_4}{a_1}(e^{a_1 \epsilon} -1), ~k= 0, \nonumber\\
{\tilde u}&=&u  + \frac{a_3}{a_1} S( e^{a_1 \epsilon} -1) +
a_4 \epsilon, ~k= 1, \label{utrosp}
\end{eqnarray}
where we assume that $a_1 \ne 0$ because the case $a_1=0$ coincides with the 
former case (\ref{str})-(\ref{utro}).
\qquad\end{proof}

We will use the
symmetry group $G_{\Delta}$ to construct invariant solutions of equation
(\ref{intur}).
To obtain the invariants of the symmetry group $G_\Delta$
 we exclude $\epsilon$ from the equations
 (\ref{str})--(\ref{utro}) or in the special case from equations (\ref{strsp})--(\ref{utrosp}).

In the first case the symmetry group $G_{\Delta}$ is very poor and we can obtain just the 
following  invariants
\begin{eqnarray}
inv_1&=&S,\label{invar1}\\
inv_2&=&u -(a_3 S +a_4)/a_2, ~~a_2 \ne 0. \nonumber
\end{eqnarray}
These invariants are useless because they do not lead to any reduction of (\ref{intur}).

In the special case (\ref{lam}) the symmetry group admits two 
functionally  independent invariants of the form
\begin{eqnarray}
inv_1&=& \log S +a t,~~a=a_1/a_2, ~a_2 \ne 0\label{invar1sp}\\
inv_2&=& u~ S^{(k-1)}.\label{invar2sp}
\end{eqnarray}

 In general the form of invariants is not unique because each function
of invariants is an invariant.
But it is possible to obtain just two non trivial functionally
independent invariants which we take in the form (\ref{invar1sp}),
(\ref{invar2sp}). 
The invariants can be used as new independent and dependent
variables in order to reduce the partial differential equation (\ref{intur})
with the special function $\lambda(S)$ defined by (\ref{lam}) to an
ordinary differential equation.

\section{The special case $\lambda(S) = \omega S^k$}\label{redsec}

Let us study a special case of equation (\ref{intur}) with
$\lambda(S)= \omega S^k, k \in R$.
The equation under investigation is now
\begin{eqnarray} \label{urav}
u_t+\frac{\sigma^2 S^2}2\frac{u_{SS}}{(1-b S^{k+1} u_{SS})^2}=0
\end{eqnarray}
with the constant $b=\rho \omega$. As usual we suggest that $\rho \in(0,1)$. 
The value of the constant $\omega$ 
depends on the corresponding option type and in our investigation it
 can be assumed that $\omega$ is an  arbitrary constant, $\omega \ne 0$. 
 The variables $S,t$ are in
the intervals
\begin{equation} \label{inter}
S>0, ~~ t \in [0, T],~~ T >0.
\end{equation}
{\rmrk The case $b=0$, i.e. $\rho=0$ or $\omega=0$ leads to the well known linear Black-Scholes model and we will exclude this case from our investigations.}

We will suppose that the denominator in equation (\ref{intur}) 
(correspondingly (\ref{urav})) is non equal to zero identically.

Let us study the denominator in the second term of the equation (\ref{urav}). It will
be equal to zero if the function $u(S,t)$ satisfies the equation
\begin{eqnarray} \label{denom}
1-b S^{k+1} u_{SS}=0.
\end{eqnarray}
The solution of this equation is a function $u_0(S,t)$
\begin{eqnarray} 
u_0(S,t)&=&\frac{1}{b k (k-1)} S^{1-k} + S c_1(t) + c_2 (t),~~b \ne 0, k\ne
0,1, \nonumber\\
u_0(S,t)&=&- \frac{1}{b} \log S + S c_1(t) + c_2 (t),~~b \ne 0, k=1,\label{nuli}\\
u_0(S,t)&=&\frac{1}{b} S \log S + S c_1(t) + c_2 (t),~~b \ne 0, k=0,\nonumber 
\end{eqnarray}
where the functions $c_1(t)$ and $c_2(t)$ are arbitrary functions of the
variable $t.$ \\[2pt]

Subsequently we will suggest that the denominator in the
second term of the equation (\ref{urav}) is  not
identically zero, i.e., a solution $u(S,t)$ is not equal to the function
$u_0(S,t)$ (\ref{nuli}) except in a discrete set of points.

Let us now introduce new invariant variables 
\begin{eqnarray}
z&=&\log S +a t, ~~ a \ne 0,\nonumber\\
v&=& u~ S^{(k-1)}. \label{newvar}
\end{eqnarray}
After this substitution equation (\ref{urav}) will be reduced to
 an ordinary differential equation
\begin{equation}
a v_z +\frac{\sigma^2}{2}\frac{v_{zz} +(1-2k)v_z -k(1-k)v}
{(1-b(v_{zz}+(1-2k)v_z -k(1-k)v))^2}=0,~~ a,b \ne 0. \label{priveq}
\end{equation}
Elementary solutions of this equation we obtain if we assume that 
$v={\rm const.}$ or $v_z={\rm const.~}$. It is easy to prove that there exists the
  trivial  solution
$v=0$ if $k \ne 0,1$, and the solutions $v={\rm const.}\ne 0$, $v={\rm const.} \ne 0$
 if $k= 0,1$ only.
The condition that the denominator in (\ref{priveq}) is non equal to zero, i.e.,
\begin{equation}
(1-b(v_{zz}+(1-2k)v_z -k(1-k)v))^2 \ne 0 \label{detnu}
\end{equation}
corresponds to equation (\ref{denom}) in new variables $ z,v $.

If the function $v(z)$ satisfies the inequality (\ref{detnu}) then we can
multiply both terms of equation (\ref{priveq}) with the denominator of the second
term.
In equation (\ref{priveq}) all coefficients are constants hence we
can reduce the order of the equation. We assume that 
$v, v_z \ne {\rm const.} $ and choose as a new independent variable $v$
and introduce as a new dependent variable $x(v)=v_z(z)$. This variable
substitution reduces equation (\ref{priveq}) to a first order differential
equation which is second order polynomial corresponding to the function
$x(v)_v$. 
Under assumption (\ref{detnu}) the set of solutions of equation (\ref{priveq})
is equivalent to a  union of solution sets 
of the following equations
\begin{eqnarray}
x&=&0~,\label{urav1}\\
x_v&=&-1 + 2\,k - \frac{\sigma^2}{4\,a\,b^2\,x^2} + \frac{1}{b\,x} +
   \frac{k (1-k)v}{x}  -
   \frac{{\sqrt{\sigma^2\,\left( \sigma^2 - 8\,a\,b\,x \right) }}}
    {4\,a\,b^2\,x^2},\nonumber \\
x_v&=&-1 + 2\,k - \frac{\sigma^2}{4\,a\,b^2\,x^2} + \frac{1}{b\,x} +
   \frac{k (1-k)v}{x} +
   \frac{{\sqrt{\sigma^2\,\left( \sigma^2 - 8\,a\,b\,x \right) }}}
    {4\,a\,b^2\,x^2} .\nonumber 
\end{eqnarray}
Equations (\ref{urav1}) are of an autonomous type if the parameter $k$ is equal to
 $k=0,1$ only. We see that these are exactly the cases in which the corresponding
 Lie-algebra (\ref{alge}) has a three dimensional  Abelian sub-algebra. 
The case $k=0$ was studied earlier in \cite{Bordag:2004}, \cite{Chmakova}. In the next 
section we will study the case $k=1$.

\section{The special case  $\lambda = \omega S$} \label{k1expl}

If we put $k=1$ in (\ref{lam}) then equation (\ref{priveq}) takes the form
\begin{equation}
v_z + q \frac{v_{zz}  -v_z }
{(1-b(v_{zz}-v_z ))^2}=0,  \label{speck1}
\end{equation}
where $q = \frac{\sigma^2}{2 a},~ a, b \ne 0. $ 
It is an autonomous equation which possesses a simple structure. We
will use this structure and introduce a more simple substitution as
described at the end of the previous section to reduce the order of equation.

One family of solutions of this equation is very easy to find. We just
suppose that the value $v_z(z)$ is equal to a constant. The equation
(\ref{speck1}) admits as a solution the value $v_z=
(-1 \pm \sqrt{q})/b$ consequently the corresponding solution
$u(s,t)$ of (\ref{urav}) with $\lambda= \omega S$ can be represented
by the formula
\begin{equation}
u(S,t)=\frac{1}{\rho \omega} \left(-1 \pm \sqrt{q}\right)(\log S + a t)+c,~ a>0, \label{constsol}
\end{equation}
where $c$ is an arbitrary constant. 

To find other families of solutions we introduce a new dependent variable
\begin{equation}
y(z)=v_z(z) \label{samy}
\end{equation} 
and assume that the denominator of the equation (\ref{speck1}) is not equal to zero,
i.e. 
\begin{equation}
v(z)\ne -\frac{z}{b}+c_1~ e^z +c_2, ~~{\rm i.e.} ~~ y(z) \ne -\frac{1}{b} +c_1~ e^z , \label{denzero}
\end{equation}
where $c_1, c_2$ are  arbitrary constants.

We multiply both terms of equation (\ref{speck1}) by the denominator of the second term and obtain
\begin{equation}
  y y_z^2 - 2 \left(y^2 +\frac{1}{b}y -\frac{q}{2 b^2} \right) y_z 
 + \left(y^2 +\frac{2}{b}y +\left(\frac{1-q}{ b^2} \right)\right) y =0,~~b \ne 0.\label{uryk1}
\end{equation} 
We denote the left hand side of this equation by $ F(y,y_z)$.
The equation (\ref{uryk1}) can possess exceptional solutions which are the solutions of a system
\begin{equation}
\frac{\partial F(y,y_z)}{\partial y_z}=0,~~ F(y,y_z)=0. \label{exepsys}
\end{equation} 
The first equation in this system defines a discriminant curve which has the form
\begin{equation}
y(z)= \frac{q}{4 b}. \label{discr}
\end{equation}
If this curve is also a solution of the original equation (\ref{uryk1}) then we
obtain an
exceptional solution.
We obtain an exceptional solution if 
$q=4$, i.e. $a=\sigma^2/8 $. It has the form 
\begin{equation}
y(z)= \frac{1}{b}. \label{exesol}
\end{equation}
This solution belongs to the family of solutions (\ref{yconst})
 by the specified value of the parameter $q$.
In all other cases the equation (\ref{uryk1})
does not possess any exceptional solutions.

Hence the set of
solutions of equation (\ref{uryk1}) is a union  of solution sets of 
following equations
\begin{eqnarray}
y&=&0,\label{ynul}\\
y&=&\left(-1 \pm \sqrt{q}\right)/b,\label{yconst}\\
y_z&=&\left( y^2 + \frac{1}{b}y   -\frac{q}{2 b^2} - \sqrt{
 \frac{\sigma^2}{2 a b^3} \left( \frac{q}{4 b}
-y\right)}\right) \frac{1}{y},~ y\ne0 \label{yzeqm}\\
y_z&=&\left( y^2 + \frac{1}{b}y   -\frac{q}{2 b^2} + \sqrt{
 \frac{q}{ b^3} \left( \frac{q}{4 b}
-y\right)}\right) \frac{1}{y}, ~ y\ne0 \label{yzeqp}
\end{eqnarray}
where one of the solutions (\ref{yconst}) is an exceptional solution (\ref{exesol}) by $q=4 $.
We denote the right hand side  of equations (\ref{yzeqm}), (\ref{yzeqp})  by
$f(y)$.  The Lipschitz condition for equations of the type $y_z=f(y)$
is satisfied in all points where the derivative $\frac{\partial
f}{\partial y}$ exists and is bounded. It is easy to see that this
condition will not be satisfied by
\begin{equation}
y=0,~~ y= \frac{q}{4 b},~~ y=\infty. \label{edin}
\end{equation}
It means that on the lines (\ref{edin}) the uniqueness of 
solutions of equations (\ref{yzeqm}), (\ref{yzeqp}) can be lost. 
We will study in detail the behavior of solutions in the neighborhood of
lines (\ref{edin}).
For this purpose we look at the equation (\ref{uryk1}) from  another point of view.
If we assume now that $z,y,y_z$ are complex variables and denote 
\begin{equation}
y(z)=\zeta, ~~ y_z(z)=w,~~ \zeta,w \in C, \label{complexsub}
\end{equation}
then the equation (\ref{uryk1}) takes the form
\begin{equation}
F(\zeta,w)= \zeta w^2 - 2 \left(\zeta^2 +\frac{1}{b}\zeta -\frac{q}{2 b^2} \right) w 
 + \left(\zeta^2 +\frac{2}{b}\zeta +\frac{1-q}{b^2}\right) \zeta =0,\label{vweq}
\end{equation}
where $ b \ne 0.$ The equation (\ref{vweq}) is an algebraic relation
in $C^2$ and defines a plane curve in this space.  
The polynomial
$F(\zeta,w)$ is an irreducible polynomial if at all roots $w_r(z)$ of $F(\zeta,w_r)$ either the
partial derivative $F_\zeta(\zeta,w_r)$ or $F_w(\zeta,w_r)$ are non equal
to zero. 
It is easy to prove that the polynomial (\ref{vweq})
is irreducible.  

We can treat equation (\ref{vweq}) as an algebraic
relation which defines a Riemann surface $\Gamma ~:~ F(\zeta,w)=0~$ of 
$w=w(\zeta)$ as a compact manifold over the $\zeta$-sphere.  
 The function $w(\zeta)$ is uniquely 
analytically extended over the Riemann surface $\Gamma$
 of two sheets over the $\zeta-$sphere.
We find all singular or branch points of $w(\zeta)$ if we study 
the roots of the first coefficient of the polynomial $F(\zeta,w)$, 
the common roots of equations
\begin{equation}
F(\zeta,w)=0,~~F_w(\zeta,w)=0,~~~ \zeta, w \in C\cup{\infty}.
\end{equation}
and the point $\zeta=\infty.$
The set of  singular or branch points consists of the points
\begin{equation}
\zeta_1=0,~~\zeta_2= \frac{q}{4 b}, ~~\zeta_3={\infty} \label{singpoint}.
\end{equation}
As expected we got the same set of points as in real case (\ref{edin}) by the study of the
Lipschitz condition but now the behavior of solutions at the points is more
visible.

The points $\zeta_2,\zeta_3$ are the branch points at which two sheets 
of $\Gamma$ are glued on. We remark that 
 \begin{equation}
w(\zeta_2)=\frac{1}{b}\left(q-4 \right) + t \frac{1}{4\sqrt{-bq}}+\cdots,~~ t^2=\zeta-\frac{q}{4 b}, \label{wz2}
\end{equation}
where $t$ is a local parameter in the neighborhood of $\zeta_2.$
For the special value of $q=4$ the value $w(\zeta_2)$ is equal to zero.

At the point $\zeta_3=\infty$ we have
$$w(\zeta)= \frac{1}{t^2} +\frac{1}{b}+t\sqrt{ \frac{-q}{4 b^3}}, ~~ t^2
=\frac{1}{\zeta},~~ \zeta \to \infty,$$
where $t$ is a local parameter in the neighborhood of $\zeta_3.$
At the point $\zeta_1=0$ the function $w(\zeta)$ has the following behavior
\begin{eqnarray}
w(\zeta)&\sim& -\frac{q}{ b^2}\frac{1}{\zeta},~\zeta \to \zeta_1=0,~ {\mbox{on the principal sheet}}, \label{pol}\\
w(\zeta) &\sim& \left(1-q\right) \zeta,~\zeta \to \zeta_1=0,~ q \ne 1,~~{\mbox{on the second sheet}}, \label{nul}\\
w(\zeta) &\sim& - 2 b^2 \zeta^2,~\zeta \to \zeta_1=0,~ q=1,~~{\mbox{on the second sheet}}. \label{nul2}
\end{eqnarray}

Any solution $w(\zeta)$ of an irreducible algebraic equation (\ref{vweq}) is
meromorphic on this compact Riemann surface $\Gamma$ of the genus 0 and has a pole of the order one
correspondingly (\ref{pol}) over the point $\zeta_1=0$ and the pole of the second
order over $\zeta_3=\infty$.
  It means also that the  meromorphic
function $w(\zeta)$ cannot be defined on a 
 manifold of less than
2 sheets over the $\zeta$ sphere.

To solve differential equations (\ref{yzeqm}) and (\ref{yzeqp})
from this point of view it is equivalent to integrate on $\Gamma$ a differential of the type $
~{\frac{{\rm d} \zeta}{w(\zeta)}}~$
and then to solve an Abel's inverse problem of degenerated type  
\begin{equation}
\int{\frac{{\rm d} \zeta}{w(\zeta)}}=z+{\rm const.}  \label{intri} 
\end{equation}
The integration can be done very easily because we can introduce a
uniformizing parameter on the Riemann surface $\Gamma$ and represent
the integral (\ref{intri}) in terms of rational functions merged 
possibly with logarithmic terms.

To realize this program we 
introduce a new variable (our uniformizing parameter $p$) in the way
\begin{eqnarray}
\zeta &=&\frac{q(1-p^2)}{4b} , \label{uniparv}\\
w&=&\frac{(1-p)(q(1+p)^2-4)}{4 b  (p+1)}.
\end{eqnarray}
Then the equations (\ref{yzeqm}) and (\ref{yzeqp}) will take the form
\begin{eqnarray}
2 q \int{\frac{p (p+1) {\rm d}p}{(p-1)(q (p+1)^2-4)}}=z+{\rm const}, \label{privedeqm}\\
2 q \int{\frac{p (p-1) {\rm d}p}{(p+1)(q (p-1)^2-4)}}=z+{\rm const}\label{privedeqp}.
\end{eqnarray}
The integration procedure of equation (\ref{privedeqm}) gives rise to the following relations
\begin{eqnarray}
&&2 q \log{(p-1)}+(q-\sqrt{q} -2) \log{\left((p+1) \sqrt{q}-2\right)} 
\label{solm1}\\
&&+(q + \sqrt{q}-2) \log{\left((p+1)\sqrt{q}+2\right)} =2 (q-1) z +c, ~~ q\ne 1,q>0\nonumber \\
&& \frac{1}{1-p} +\frac{1}{4}\log{\frac{(p+3)^3}{(p-1)^5}}=z+ c, ~~ q=1.\label{solm2}\\
&&2 \sqrt{(-q)} \arctan{((p+1)\sqrt{(-q)}/2)} -2 q \log{(p-1)}\nonumber \\
&&+(2-q)\log{(4 - q (p+1)^2)}=2 (1-q) z +c, ~~ q<0, \label{solm3}
\end{eqnarray}
where $c$ is an arbitrary constant.
The equation (\ref{privedeqp}) leads to
\begin{eqnarray}
&&2 q \log{(p+1)}+(q+\sqrt{q} -2) \log{\left((p-1) \sqrt{q}-2\right)} 
\label{solp1}\\
&&+(q - \sqrt{q}-2) \log{\left((p-1)\sqrt{q}+2\right)} =2 (q-1) z +c, ~~ q\ne 1,q>0\nonumber \\
&& \frac{1}{p+1} +\frac{1}{4}\log{\frac{(p-3)^3}{(p+1)^5}}=z+ c, ~~ q=1.\label{solp2}\\
&&-2 \sqrt{(-q)} \arctan{((p-1)\sqrt{(-q)}/2)} -2 q \log{(1+p)}\nonumber \\
&&+(2-q)\log{(4 - q (p-1)^2)}=2 (1-q) z +c, ~~ q<0. \label{solp3}
\end{eqnarray}
where $c$ is an arbitrary constant.

The relations (\ref{solm1})-(\ref{solp3}) are first order ordinary differential equations
because of the 
 substitutions (\ref{complexsub}) and (\ref{samy}) we have 
\begin{equation}
p=\sqrt{ 1 -\frac{4 b}{q }  v_z }. \label{vzp}
\end{equation}
All these results can be collected to the following  theorem.

\begin{theorem}
The equation (\ref{speck1}) for arbitrary values of the
parameters $q, b\ne 0$
can be reduced to the set of first order differential equations which consists
of the equations 
\begin{equation}
 v_z=0,~~v_z=(-1 \pm \sqrt{q})/b \label{trivisol}
\end{equation}
and equations (\ref{solm1})-(\ref{solp3}). The complete set of solutions of
the equation (\ref{speck1}) coincides with the union of solutions of these equations.
\end{theorem}

To solve equations (\ref{solm1})-(\ref{solp3}) 
exactly we should first invert these formulas 
in order to obtain an
exact representation $p$ as a function of $z$.
If an exact formula for the function $p=p(z)$  is found we can use the 
substitution (\ref{vzp}) to obtain an explicit ordinary differential equation of
the type $v_z(z)=f(z)$ or another suitable type
and if it possible then to integrate the final equation. 

But even on the first step we would not be able to do this for an
 arbitrary value of the parameter $q$.  It means we have just implicit
 representations for the solutions of the equation (\ref{speck1}) as
 solutions of the implicit first order differential equations
 (\ref{solm1})-(\ref{solp3}).

\subsection{Exact invariant solutions in case of 
 a fixed relation between variables $S$ and $t$}

For a special value of the parameter $q$ we can invert the equations (\ref{solm1}) 
and (\ref{solp1}). Let us take $q=4$, i.e., the relation between variables $S,t$ is
fixed in the form
\begin{equation} 
z=\log S +\frac{\sigma^2}{8}t.  \label{podszotst}
\end{equation}
In this case the equation (\ref{solm1})
takes the form
\begin{equation}
(p-1)^2(p+2)=c~ \exp{(3 z/2)} \label{cubm}
\end{equation}
and correspondingly the equation (\ref{solp1}) the form
\begin{equation}
(p+1)^2(p-2)=c~ \exp{(3 z/2)}, \label{cubp}
\end{equation}
where $c$ is an arbitrary constant.
It is easy to see that the equations (\ref{cubm}) and (\ref{cubp}) are connected by
a transformation
\begin{equation}
p\to -p, ~~ c\to -c. \label{involut} 
\end{equation}
This  symmetry arises from the symmetry of the underlining  Riemann surface  
$\Gamma$ (\ref{vweq})
and corresponds to a change of the sheets on  $\Gamma$.

\begin{theorem}
The second order differential equation
\begin{equation}
v_z + 4 \frac{v_{zz} -v_z }
{(1-b(v_{zz}-v_z ))^2}=0,  \label{speck1q4}
\end{equation}
is exactly integrable for an arbitrary value of the parameter $b$. The
 complete set of solutions
for $b \ne 0$ is given by the union of
solutions
(\ref{solvzplus}), (\ref{resh1}) -(\ref{resh32}) and solutions 
\begin{equation}
v(z)=d,~~v(z)=-\frac{3}{b} z+d,~~v(z)=\frac{1}{b} z+d,\label{trivisol1}
\end{equation}
where $d$ is an arbitrary constant. The last solution in (\ref{trivisol1}) corresponds to the
 exceptional solution of equation (\ref{uryk1}).\\
For $b=0$  equation (\ref{speck1q4}) is linear and its solutions are given by $v(z)=d_1 +
d_2 \exp{(3 z/4 )}$, where $d_1,d_2$ are arbitrary constants.
\end{theorem} 
 
\begin{figure}[ht]
\vspace{0.5cm}
\begin{minipage}[t]{12cm}
\includegraphics[width=12cm]{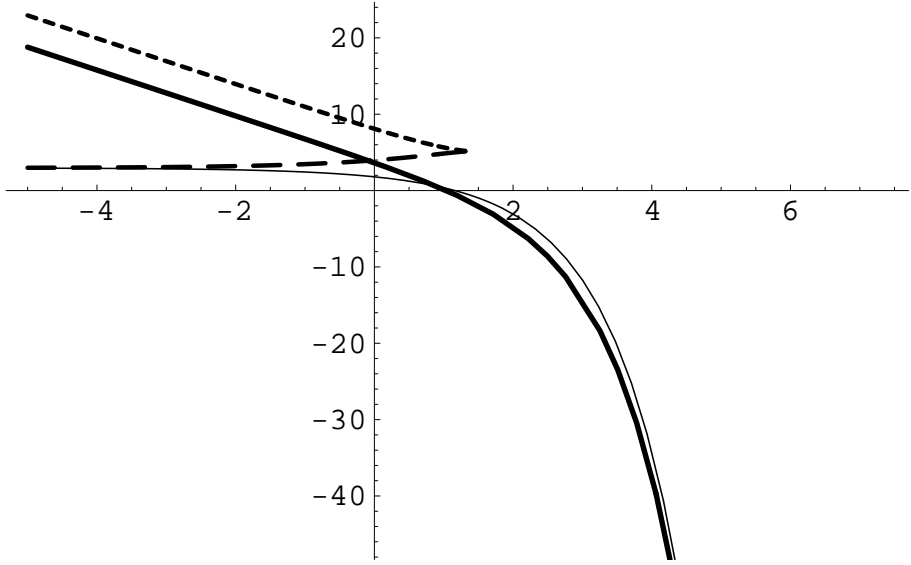}
\caption{\label{solfi2} Plot of the solution $v(z)$, (\ref{solvzplus}),
(thick solid line),
$v_1(z)$, (\ref{resh1}), (short dashed line), 
$v_2(z)$,  (\ref{resh2}), (long dashed line) and the third solution 
$v_{3,1}(z),v_{3,2}(z) $,  (\ref{resh31}), (\ref{resh32}),
which is represented by the thin solid line.
The parameters takes the values $~|c|=1,q=4, d=0, b=1$ and the variable $z\in (-5,4.5)$. }
\end{minipage} 
\end{figure}

\begin{proof}
Because of the symmetry (\ref{involut}) it is sufficient to study eithr the
equations (\ref{cubm}) or (\ref{cubp}) for $c\in R$ or both these equations 
for $c>0$. The value $c=0$ can
be excluded because it complies with the constant value of $p(z) $ and
correspondingly constant value of $v_z(z)$, but all such cases are
studied before and the solutions are given by (\ref{trivisol1}).

We will study equation (\ref{cubp}) in case $c \in R \setminus \{0\}$
and obtain on this way the complete class of exact solutions for equations
(\ref{cubm})-(\ref{cubp}).

Equation (\ref{cubp}) for $c>0$ has a one real root only. It leads to
an ordinary differential equation of the form
\begin{eqnarray}
 v_z(z)&=&- \frac{1}{b} -
   \frac{2^{\frac{2}{3}}}{b \,
      {\left( 2 + c\,e^{\frac{3\,z}{2}} +
          {\sqrt{4\,c\,e^{\frac{3\,z}{2}} + c^2\,e^{3\,z}}} \right) }^
       {\frac{2}{3}}} \nonumber\\
&&- \frac{
      {\left( 2 + c\,e^{\frac{3\,z}{2}} +
          {\sqrt{4\,c\,e^{\frac{3\,z}{2}} + c^2\,e^{3\,z}}} \right) }^
       {\frac{2}{3}}}{b \,2^{\frac{2}{3}}}, ~~c>0. \label{reshv}
\end{eqnarray}

Equation (\ref{reshv}) can be exactly integrated if we use an  Euler substitution
and introduce a new independent variable
\begin{equation}
\tau=2 + c\,e^{\frac{3\,z}{2}} + {\sqrt{4\,c\,e^{\frac{3\,z}{2}} +
       c^2\,e^{3\,z}}}. \label{podstau1}
\end{equation}

The corresponding solution is given by
\begin{eqnarray}
v(z) &=& - \frac{2^{\frac{2}{3}}}
     {{b \left( 2 + c\,e^{\frac{3\,z}{2}} +
           {\sqrt{4\,c\,e^{\frac{3\,z}{2}} + c^2\,e^{3\,z}}} \right) }^
        {\frac{2}{3}}} 
- \frac{{\left( 2 + c\,e^{\frac{3\,z}{2}} +
         {\sqrt{4\,c\,e^{\frac{3\,z}{2}} + c^2\,e^{3\,z}}} \right) }^
      {\frac{2}{3}}}{b 2^{\frac{2}{3}}} \nonumber\\
&-& \frac{2}{b} \log \left( \frac{2^{\frac{1}{3}}}
      {{\left( 2 + c\,e^{\frac{3\,z}{2}} +
           {\sqrt{4\,c\,e^{\frac{3\,z}{2}} + c^2\,e^{3\,z}}} \right) }^{\frac{1}{3}}}  
+ \frac{{\left( 2 + c\,e^{\frac{3\,z!
 }{2}} +
         {\sqrt{4\,c\,e^{\frac{3\,z}{2}} + c^2\,e^{3\,z}}} \right) }^
      {\frac{1}{3}}}{2^{\frac{1}{3}}} -2 \right) \nonumber\\ 
&+&d,\label{solvzplus}
\end{eqnarray}
where $d \in R$ is an arbitrary constant.\\
If in the right hand side of equation (\ref{cubp}) the parameter $c $ 
satisfies the inequality $c<0$ and the variable $z$ chosen in the region
\begin{equation}
 z \in \left(-\infty,\frac{4}{3}\ln{\frac{2}{|c|}} \right)
\end{equation}
then the equation on $p$  possesses maximal three real roots.

These three roots of cubic equation  (\ref{cubp}) give rise to
 three differential equations of the type
$v_z=(1- p^2(z))/b$. The equations can be exactly solved and we find
correspondingly three solutions $v_i(z),~ i=1,2,3$.

The first solution is given by the expression
\begin{eqnarray}
v_1(z)  &=&\frac{z}{b}-\frac{2}{b}\,\cos \left(\frac{2}{3}\,\arccos 
\left(1 - \frac{|c|}{2}\,e^{\frac{3\,z}{2}}\right)\right) \label{resh1}\\
&-& \frac{4}{3 b}\,\log \left(1 + 2\,\cos 
\left(\frac{1}{3}\,\arccos 
\left(1 - \frac{|c|}{2}\,e^{\frac{3\,z}{2}}\right)\right)\right)\nonumber\\
&-& \frac{16}{3 b}
     \,\log \left(\sin 
\left(\frac{1}{6}\,\arccos 
\left(1 - \frac{|c|}{2}\,e^{\frac{3\,z}{2}}\right)\right)\right)+d,\nonumber
\end{eqnarray}
where $d \in R$ is an arbitrary constant.
The second solution is given by the formula
\begin{eqnarray}
v_2(z) &=&\frac{z}{b}-\frac{2}{b}\,\cos \left(\frac{2}{3} \pi +\frac{2}{3}\,\arccos 
\left(-1 + \frac{|c|}{2}\,e^{\frac{3\,z}{2}}\right)\right) \label{resh2}\\
&-& \frac{4}{3 b}\,\log \left(1 + 2\,\cos 
\left(\frac{1}{3} \pi +\frac{1}{3}\,\arccos 
\left(-1 + \frac{|c|}{2}\,e^{\frac{3\,z}{2}}\right)\right)\right)\nonumber\\
&-& \frac{16}{3 b}
     \,\log \left(\sin 
\left(\frac{1}{6} \pi +\frac{1}{6}\,\arccos 
\left(-1 + \frac{|c|}{2}\,e^{\frac{3\,z}{2}}\right)\right)\right)+d,\nonumber
\end{eqnarray}
where $d \in R$ is an arbitrary constant.
The first and second solutions are defined  up to the point 
$z=\frac{4}{3}\ln{\frac{2}{|c|}}$ where they coincide (see Fig.~\ref{solfi2}).

The third solution for  $z< \frac{4}{3}\ln{\frac{2}{|c|}}$ is given by the formula
\begin{eqnarray}
v_{3,1}(z) &=&\frac{z}{b} -\frac{2}{b}\,\cos \left(\frac{2}{3}\,\arccos 
\left(-1 + \frac{|c|}{2}\,e^{\frac{3\,z}{2}}\right)\right) \label{resh31}\\
&-& \frac{4}{3 b}\,\log \left(-1 + 2\,\cos 
\left(\frac{1}{3}\,\arccos 
\left(-1 + \frac{|c|}{2}\,e^{\frac{3\,z}{2}}\right)\right)\right)\nonumber\\
&-& \frac{16}{3 b}
     \,\log \left(\cos 
\left(\frac{1}{6}\,\arccos 
\left(-1 + \frac{|c|}{2}\,e^{\frac{3\,z}{2}}\right)\right)\right)+d,\nonumber
\end{eqnarray}
where $d \in R$ is an arbitrary constant.
In case $z> \frac{4}{3}\ln{\frac{2}{|c|}}$ the polynomial (\ref{cubp})
has a one real root and the corresponding solution can be represented by the formula
\begin{eqnarray}
v_{3,2}(z)&=& \frac{z}{b} - \frac{2}{b}\,\cosh \left(\frac{2}{3}
          {\rm arccosh}\left(-1 +
          \frac{|c|}{2}\,e^{\frac{3\,z}{2}}\right)\right)\label{resh32}\\
&-&
    \frac{16}{3\,b}\,\log \left(\cosh \left(\frac{1}{6}{\rm arccosh}\left(-1 +
            \frac{|c|}{2}\,e^{\frac{3\,z}{2}}\right)\right)\right)\nonumber \\
&-&
   \frac{4}{3\,b}\,\log \left(-1 + 2\,\cosh \left(\frac{1} {3} {\rm arccosh} \left(-1 +
              \frac{|c|}{2}\,e^{\frac{3\,z}{2}}\right)\right)\right) +d.\nonumber
\end{eqnarray}

The third solution is represented by formulas $v_{3,2}(z)$ and  $v_{3,1}(z)$
for different values of the variable $z$.
\qquad\end{proof} 

One of the sets of solutions (\ref{solvzplus}), (\ref{resh1})
-(\ref{resh32}) for fixed parameters $b,c,d$ is represented in
Fig.~\ref{solfi2}. The first solution (\ref{reshv}) and the third
solution given by both (\ref{resh31}) and (\ref{resh32}) are defined for any
values of $z$. The solutions $v_1(z)$ and $v_2(z)$ cannot be
continued after the point
$z=\frac{4}{3}\ln{\frac{2}{|c|}}$ where they coincide.

\begin{figure}[ht]
\vspace{0.5cm}
\begin{minipage}[t]{12cm}
\includegraphics[width=12cm]{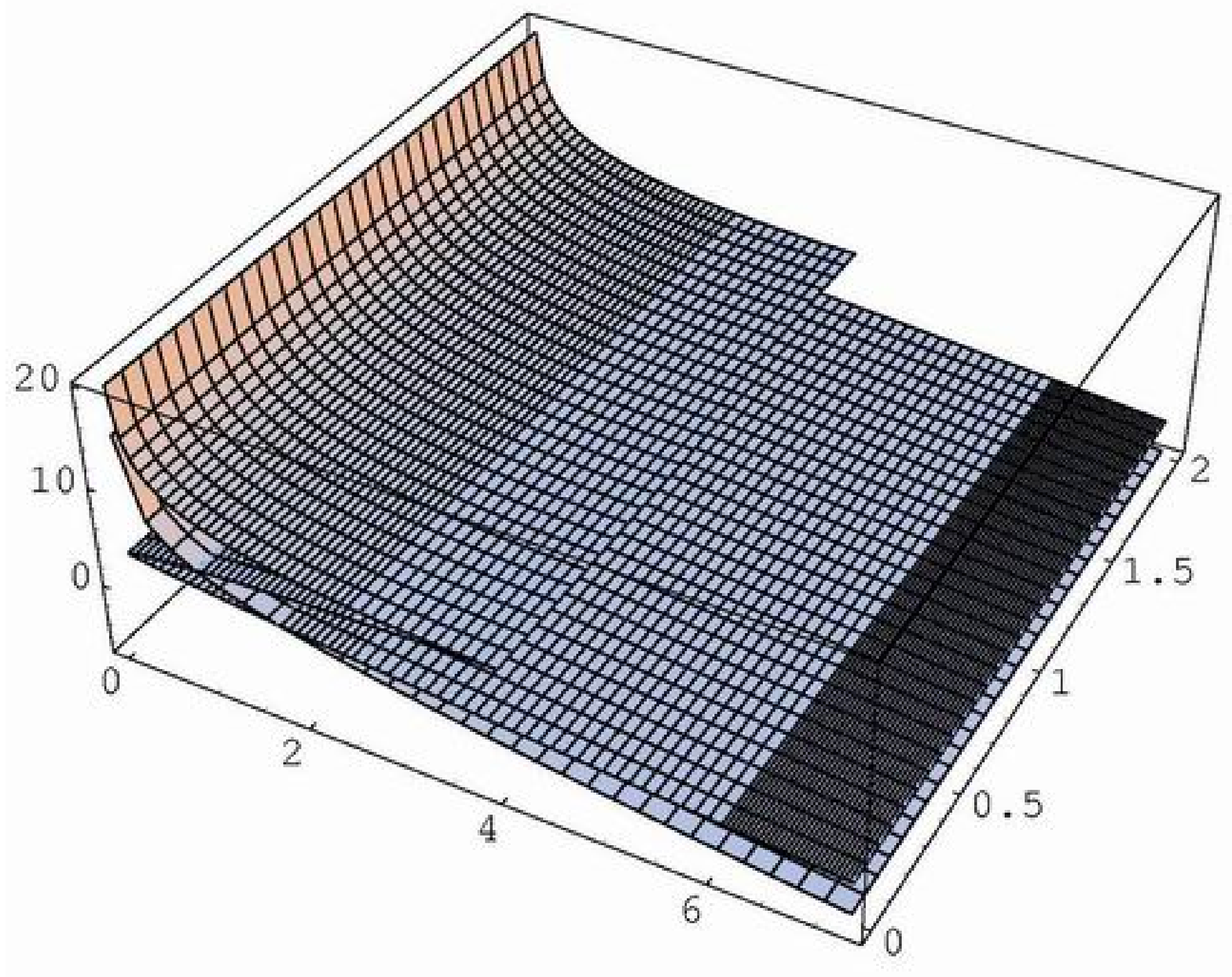}
\caption{\label{solfi3} Plot of solutions $u(S,t),u_1(S,t),u_2(S,t),u_{3,1}(S,t),u_{3,2}(S,t) $  for the
  parameters $|c|=0.5,q=4,b=1.0, d=0$. The variables $S,t$ lie in intervals $S\in (0,9)$ and $t\in [0,2.0]$. 
All invariant solutions change slowly in $t$-direction.}
\end{minipage} 
\end{figure}

If we put in mind that $z=\log S +\frac{\sigma^2}{8}t$ and $v(z)=u(S,t)$  
we can represent exact
invariant solution of equation (\ref{urav}).
The solution (\ref{solvzplus}) gives rise to an invariant solution
$u(S,t)$ in the form
\begin{eqnarray}
u(S,t) = - 
      \frac{1}{\omega \rho}\left( 1 + c\,S^{\frac{3}{2}} e^{\frac{3 \sigma^2}{16} t} +
           {\sqrt{2\,c\,S^{\frac{3}{2}} e^{\frac{3 \sigma^2}{16}t } + c^2\,S^3 e^{\frac{3 \sigma^2}{8}t}}} 
\right)^{-\frac{2}{3}}   \nonumber\\
- \frac{1}{\omega \rho} {\left( 1 + c\,S^{\frac{3}{2}} e^{\frac{3 \sigma^2}{16}t} +
         {\sqrt{2\,c\,S^{\frac{3}{2}} e^{\frac{3 \sigma^2}{16}t} + c^2\,S^3 e^{\frac{3 \sigma^2}{8}t}}} \right) }^
      {\frac{2}{3}}~~~~~\label{solustl}\\
-  \frac{2}{\omega \rho}\log \left( 
      {\left( 1 + c\,S^{\frac{3}{2}} e^{\frac{3 \sigma^2}{16}t} +
           {\sqrt{2\,c\,S^{\frac{3}{2}} e^{\frac{3 \sigma^2}{16}t} + c^2\,S^3 e^{\frac{3 \sigma^2}{8}t}}} \right) }^
        {-\frac{1}{3}}  \right. \nonumber\\
\left.
~~~~~~~~~~~~~~~~~~+{\left( 1 + c\,S^{\frac{3}{2}} e^{\frac{3 \sigma^2}{16}t} +
         {\sqrt{2\,c\,S^{\frac{3}{2}} e^{\frac{3 \sigma^2}{16}t} + c^2\,S^3 e^{\frac{3 \sigma^2}{8}t}}} \right) }^
      {\frac{1}{3}} -2 \right)  
+d  \nonumber
\end{eqnarray}
where $d \in R$,  $c>0$.

In case $c<0$ we can obtain correspondingly three solutions if
\begin{equation} \label{ogranst}
 0<S\le\left(\frac{2}{|c|}\right)^{4/3} \exp{\left(-\frac{\sigma^2}{8}  t\right)}.\end{equation}
The first solution is represented by
\begin{eqnarray}
u_1(S,t)  &=&\frac{1}{\omega \rho}\left(\log {S} + \frac{\sigma^2}{8}  t\right)-\frac{2}{\omega \rho}\,\cos \left(\frac{2}{3}\,\arccos 
\left(1 - \frac{|c|}{2}\,
S^{\frac{3}{2}} e^{\frac{3 \sigma^2}{16}t}
\right)\right) \nonumber\\
&-& \frac{4}{3 \omega \rho}\,\log \left(1 + 2\,\cos 
\left(\frac{1}{3}\,\arccos 
\left(1 - \frac{|c|}{2}\,
S^{\frac{3}{2}} e^{\frac{3 \sigma^2}{16}t}
\right)\right)\right)\label{solustl1}\\
&-& \frac{16}{3 \omega \rho}
     \,\log \left(\sin 
\left(\frac{1}{6}\,\arccos 
\left(1 - \frac{|c|}{2}\,
S^{\frac{3}{2}} e^{\frac{3 \sigma^2}{16}t}
\right)\right)\right)+d,\nonumber
\end{eqnarray}
where $d \in R$,  $c<0$.
The second solution is given by the formula
\begin{eqnarray}
u_2(S,t) &=&\frac{1}{\omega \rho}\left(\log {S} + \frac{\sigma^2}{8}  t\right)
-\frac{2}{\omega \rho}\,\cos \left(\frac{2}{3} \pi +\frac{2}{3}\,\arccos 
\left(-1 + \frac{|c|}{2}\,
S^{\frac{3}{2}} e^{\frac{3 \sigma^2}{16}t}
\right)\right) \nonumber\\
&-& \frac{4}{3 \omega \rho}\,\log \left(1 + 2\,\cos 
\left(\frac{1}{3} \pi +\frac{1}{3}\,\arccos 
\left(-1 + \frac{|c|}{2}\,
S^{\frac{3}{2}} e^{\frac{3 \sigma^2}{16}t}\right)\right)\right) \label{solustl2}\\
&-& \frac{16}{3 \omega \rho}
     \,\log \left(\sin 
\left(\frac{1}{6} \pi +\frac{1}{6}\,\arccos 
\left(-1 + \frac{|c|}{2}\,
S^{\frac{3}{2}} e^{\frac{3 \sigma^2}{16}t}\right)\right)\right)+d. \nonumber
\end{eqnarray}
where $d \in R$,  $c<0$.
The first and second solutions are defined for the variables under conditions
(\ref{ogranst}).
They coincide along the curve 
$$S=\left(\frac{2}{|c|}\right)^{4/3} \exp{\left(-\frac{\sigma^2}{8}  t \right)}$$
and cannot be continued further.

The third solution is defined by 
\begin{eqnarray}
u_{3,1}(S,t) &=&\frac{1}{\omega \rho}\left(\log {S} + \frac{\sigma^2}{8}  t\right)
-\frac{2}{\omega \rho}\,\cos \left(\frac{2}{3}\,\arccos 
\left(-1 + \frac{|c|}{2}\,
S^{\frac{3}{2}} e^{\frac{3 \sigma^2}{16}t}
\right)\right) \nonumber\\
&-& \frac{4}{3 \omega \rho}\,\log \left(-1 + 2\,\cos 
\left(\frac{1}{3}\,\arccos 
\left(-1
  + \frac{|c|}{2}\,
S^{\frac{3}{2}} e^{\frac{3 \sigma^2}{16}t}
 \right)\right)\right)\label{solustl31}\\
&-& \frac{16}{3 \omega \rho}
     \,\log \left(\cos \left(\frac{1}{6}\,\arccos \left(-1 + \frac{|c|}{2}\,
S^{\frac{3}{2}} e^{\frac{3 \sigma^2}{16}t}
\right)\right)\right)+d, \nonumber
\end{eqnarray}
where $d \in R$ and $S,t$ satisfied the condition (\ref{ogranst}).

In case $\log {S} + \frac{\sigma^2}{8}  t> \frac{4}{3}\ln{\frac{2}{|c|}}$ the
third solution can be represented by the formula
\begin{eqnarray}
u_{3,2}(S,t)&=& \frac{1}{\omega \rho}\left(\log {S} + \frac{\sigma^2}{8}  t\right)
- \frac{2}{\omega \rho}\,\cosh \left(\frac{2}{3}
          {\rm arccosh}\left(-1 +
          \frac{|c|}{2}\,
S^{\frac{3}{2}} e^{\frac{3 \sigma^2}{16}t}
\right)\right)\nonumber \\
&-&
    \frac{16}{3\,\omega \rho}\,\log \left(\cosh \left(\frac{1}{6}{\rm arccosh}\left(-1 +
            \frac{|c|}{2}\,
S^{\frac{3}{2}} e^{\frac{3 \sigma^2}{16}t}
\right)\right)\right) \label{solustl32} \\
&-&
   \frac{4}{3\,\omega \rho}\,\log \left(-1 + 2\,\cosh \left(\frac{1} {3} {\rm arccosh} \left(-1 + \frac{|c|}{2}\,
S^{\frac{3} {2}} e^{\frac{3 \sigma^2}{16}t}
\right)\right)\right) +d.\nonumber
\end{eqnarray}

The solution $u(S,t)$ (\ref{solustl}) and the third solution given by $u_{3,1}$,
$u_{3,2}$ (\ref{solustl31}),(\ref{solustl32}) are defined for all values of
variables $t$ and $S>0$. They have a common intersection curve of the type $S={\rm const.}~ \exp(-{\frac{ \sigma^2}{8}t} )$. The typical behavior of all these invariant solutions is represented in Fig.~\ref{solfi3}.

Previous results can be summed up in the following theorem describing the set
of invariant solutions of equation (\ref{intur}).

\begin{theorem}
\begin{remunerate}
\item The equation (\ref{intur}) possesses invariant solutions for the special
  form of the function $\lambda(S)$ given by (\ref{lam}) only. 
\item In case (\ref{lam}) the invariant
 solutions of  equation (\ref{intur}) are defined by 
ordinary differential equations (\ref{urav1}).
In special cases $k=0,1$  equations  (\ref{urav1}) are of an autonomous type.
\item If $\lambda (S)= \omega S $, i.e. $k=1$, then  the invariant solutions of equation (\ref{intur}) can be defined by 
the set of first order ordinary differential equations (\ref{solm1})--(\ref{solp3}) and
  equation (\ref{trivisol}). \\
If additionally the parameter $q=4$, or equivalent in the first invariant 
(\ref{invar1sp}) we chose $a=\sigma^2/8$ then the complete set of invariant solutions
  (\ref{intur}) can
  be found exactly. This set of invariant solutions is given by formulas
  (\ref{solustl})--(\ref{solustl32}) and by solutions
$$u(S,t)=d,~~u(S,t)=-3/b ~(\log{S}+\sigma^2 t/8),~~u(S,t)=1/b ~(\log{S}+\sigma^2
  t/8),$$
where $d$ is an arbitrary constant.
This set of invariant solutions is unique up to the transformations of the symmetry group $G_{\Delta}$
  given by theorem \ref{symteor}.
\end{remunerate}
\end{theorem}

The solutions $u(S,t)$ (\ref{solustl}), $u_1(S,t)$ (\ref{solustl1}),
$u_2(S,t)$ (\ref{solustl2}), $u_{3,1}(S,t)$ (\ref{solustl31}),
$u_{3,2}(S,t)$ (\ref{solustl32}), have no one counterpart in a linear
case. If the parameter $\rho \to 0$ then equation (\ref{intur}) and
correspondingly equation (\ref{urav}) will be reduced to the linear
Black-Scholes equation but solutions (\ref{solustl})-(\ref{solustl32})
which we obtained here will be completely blown up by $\rho \to 0$
because of the factor $1/b=1/(\omega \rho)$ in the formulas
(\ref{solustl})-(\ref{solustl32}).  This phenomena was described as
well in \cite{Bordag:2004}, \cite{Chmakova} for the invariant
solutions of equation (\ref{urav}) with $k=0$.

\nocite{Chmakova}
\nocite{Bordag:2004}
\nocite{Frey:illicuidity}
\nocite{FreyPatie} 
\nocite{Frey:perfect} 
\nocite{SchonbucherWilmotta}
\nocite{Gaeta} 
\nocite{Stephani} 
\nocite{Olver} 
\nocite{Ibragimov}
\nocite{Ovsiannikov}
\nocite{Lie}
\bibliographystyle{siam}
\bibliography{bordagchmakova}
\end{document}